\input amstex
\magnification=1200
\documentstyle{amsppt}
\NoRunningHeads
\NoBlackBoxes
\define\Mat{\operatorname{Mat}}
\topmatter
\title Tactics, dialectics, representation theory
\endtitle
\author Denis V. Juriev
\endauthor
\affil ul.Miklukho-Maklaya 20-180, Moscow 117437 Russia\linebreak
(e-mail: denis\@juriev.msk.ru)\linebreak
\ \linebreak
\ \linebreak
\endaffil
\date math.HO/0001032\newline January 06, 2000\enddate
\abstract\nofrills\linebreak
\ \linebreak
\ \linebreak

This article is devoted to the tactical game theoretical
interpretation of dialectics. Dialectical games are considered as abstractly
as well as models of the internal dialogue and reflection. The models
related to the representation theory (representative dynamics) are specially 
investigated in detail, they correlate with the hypothesis on the dialectical 
features of human thinking in general and mathematical thought (the
constructing of a solution of mathematical problem) in particular.

\ \linebreak
\ \linebreak
\endabstract
\toc
\specialhead Introduction\endspecialhead
\head I. \ \ Tactics\endhead
\subhead 1.1. Interactive games\endsubhead
\subhead 1.2. Dialogues and verbalizable games\endsubhead
\subhead 1.3. Tactical games\endsubhead
\head II. \ \ Dialectics and dialectical games\endhead
\head III. \ \ Dialectical games and representation theory\endhead
\subhead 1.1. Representative dynamics\endsubhead
\subhead 1.2. Tactical representative dynamics as dialectical game\endsubhead
\specialhead Conclusions\endspecialhead
\endtoc
\endtopmatter
\document
\newpage
\head Introduction\endhead

The mathematical formalism of interactive games, which extends one of
ordinary games (see e.g.[1]) and is based on the concept of an interactive 
control, was recently proposed by the author [2] to take into account 
the complex composition of controls of a real human person, which are often 
complicated couplings of his/her cognitive and known controls with 
the unknown subconscious behavioral reactions. In the article [3] the
dialogues as psycholinguistic phenomena were described in the interactive
game theoretical terms. It allowed to transfer from the pure interactivity to
the essentially new and much complex phenomenon of tactics [4]. The tactical
games were introduced together with some useful methods of their constructing,
examples were considered and applications specified. Generally, tactics may
be regarded as "an art to manipulate the unknown, which is manifested by
the interactivity, without making it known" [4], in particular, tactical 
actions are not completely understood even by the acting person.

The present article is devoted to the abstract foundations of tactics, 
precisely, to the interactive game theoretical interpretation of dialectics. 
The least is defined as a logical self-describing tactical game. A proper 
subclass of tactical games is extracted, they are the dialectical games,
the tactical games, which are tactical extensions of dialectics. And if 
tactics is thought as an art to manipulate the unknown without making it 
known, dialectics may be thought as an art to comprehend such manipulations. 

An interesting type of dialectical games is related to the representation
theory, with a strong accent on its inverse problems and representative
dynamics. The special attention is paid to such games as well as to the 
various sides of interaction of tactics, dialectics and representation 
theory that is the main theme of this article. 

The present article finishes two-year researches of the author on the
mathematical formulation of the game theoretical background for interactivity 
and tactics. It seems that the picture is complete now and all main concepts
are introduced. However, many technical details should be clarified but
the practice is essential here. Certainly, only the practice may determine
the directions of future theoretical studies. Mathematical proof of any 
nontrivial statement claims a lot of time and efforts, and when the initial 
direction was misguiding the whole activity is simply senseless. So let's 
``learn to labor and to wait''.

\head I. Tactics\endhead

\subhead 1.1. Interactive games [2]\endsubhead

\definition{Definition 1 [2]} An {\it interactive system\/} (with $n$
{\it interactive controls\/}) is a control system with $n$ independent 
controls coupled with unknown or incompletely known feedbacks (the feedbacks
as well as their couplings with controls are of a so complicated nature that 
their can not be described completely). An {\it interactive game\/} is a game 
with interactive controls of each player.
\enddefinition

Below we shall consider only deterministic and differential interactive
systems. In this case the general interactive system may be written in the 
form:
$$\dot\varphi=\Phi(\varphi,u_1,u_2,\ldots,u_n),\tag1$$
where $\varphi$ characterizes the state of the system and $u_i$ are
the interactive controls:
$$u_i(t)=u_i(u_i^\circ(t),\left.[\varphi(\tau)]\right|_{\tau\leqslant t}),$$
i.e. the independent controls $u_i^\circ(t)$ coupled with the feedbacks on
$\left.[\varphi(\tau)]\right|_{\tau\leqslant t}$. One may suppose that the
feedbacks are integrodifferential on $t$.

However, it is reasonable to consider the {\it differential interactive
games}, whose feedbacks are purely differential. It means that
$$u_i(t)=u_i(u_i^\circ(t),\varphi(t),\ldots,\varphi^{(k)}(t)).$$
A reduction of general interactive games to the differential ones via the 
introducing of the so-called {\it intention fields\/} was described in [2]. 
Below we shall consider the differential interactive games only if the opposite 
is not specified explicitely.

The interactive games introduced above may be generalized in the following 
ways. 

The first way, which leads to the {\it indeterminate interactive games},
is based on the idea that the pure controls $u_i^\circ(t)$ and the 
interactive controls $u_i(t)$ should not be obligatory related in the
considered way. More generally one should only postulate that there are
some time-independent quantities $F_\alpha(u_i(t),u_i^\circ(t),\varphi(t),
\ldots,\varphi^{(k)}(t))$ for the independent magnitudes $u_i(t)$ and 
$u_i^\circ(t)$. Such claim is evidently weaker than one of Def.1. For 
instance, one may consider the inverse dependence of the pure and 
interactive controls: $u_i^\circ(t)=u_i^\circ(u_i(t),\varphi(t),\ldots,
\varphi^{(k)}(t))$.

The second way, which leads to the {\it coalition interactive games}, is
based on the idea to consider the games with coalitions of actions and to
claim that the interactive controls belong to such coalitions. In this case
the evolution equations have the form
$$\dot\varphi=\Phi(\varphi,v_1,\ldots,v_m),\tag2$$
where $v_i$ is the interactive control of the $i$-th coalition. If the 
$i$-th coalition is defined by the subset $I_i$ of all players then
$$v_i=v_i(\varphi(t),\ldots,\varphi^{(k)}(t),u^\circ_j| j\in I_i).$$
Certainly, the intersections of different sets $I_i$ may be non-empty so
that any player may belong to several coalitions of actions. Def.1 gives the
particular case when $I_i=\{i\}$.

\remark{Remark 1} One is able to consider interactive games of discrete time 
in the similar manner.
\endremark

\remark{Remark 2} If one suspect that the explicit dependence of the feedbacks
on the derivatives of $\varphi$ is not correct because they are determined
via the evolution equations governed by the interactive controls, it is
reasonable to use the inverse dependence of pure and interactive controls.
\endremark

Interactive games are games with incomplete information by their nature.
However, this incompleteness is in the unknown feedbacks, not in the 
unknown states. The least situation is quite familiar to specialists in
game theory and there is a lot of methods to have deal with it. For
instance, the unknown states are interpreted as independent controls of
the virtual players and some muppositions on their strategies are done.
To transform interactive games into the games with an incomplete information
on the states one can use the following trick, which is called the 
$\varepsilon$-representation of the interactive game.

\definition{Definition 2} The $\varepsilon$-representation of the 
differential interactive game is a representation of the interactive controls 
$u_i(t)$ in the form
$$u_i(t)=u_i(u^\circ_i(t),\varphi(t),\ldots\varphi^{(k)}(t);\varepsilon_i(t))$$
with the {\sl known\/} function $u_i$ of its arguments $u_i^\circ$,
$\varphi,\ldots,\varphi^{(k)}$ and $\varepsilon_i$, whereas 
$$\varepsilon_i(t)=\varepsilon_i(u^\circ_i(t),\varphi(t),\ldots,\varphi^{(k)}
(t)$$ 
is the {\sl unknown\/} function of $u_i^\circ$ and $\varphi,\ldots,
\varphi^{(k)}$.
\enddefinition

$\varepsilon_i$ are interpreted as parameters of feedbacks and, thus, 
characterize the internal {\sl existential\/} states of players. It motivates 
the notation $\varepsilon$. Certainly, $\varepsilon$-parameters are not
really states being the unknown functions of the states and pure controls,
however, one may sometimes to apply the standard procedures of the theory
of games with incomplete information on the states. For instance, it is
possible to regard $\varepsilon_i$ as controls of the virtual players.
The na{\"\i}vely introduced virtual players only double the ensemble of 
the real ones in the interactive games but in the coalition interactive 
games the collective virtual players are observed. More sophisticated
procedures generate ensembles of virtual players of diverse structure.

Precisely, if the derivatives of $\varphi$ are excluded from the feedbacks (at 
least, from the interactive controls $u_i$ as functions of the pure controls,
states and the $\varepsilon$-parameters) the evolution equation will have 
the form
$$\dot\varphi(t)=\Phi(\varphi,u_1(u^\circ_1(t),\varphi(t);\varepsilon_1(t)),
\ldots,u_n(u_n^\circ(t),\varphi(t);\varepsilon_n(t))),\tag3$$
so it is consistent to regard the equations as ones of the controlled
system with the ordinary controls $u_1,\ldots,u_n,\varepsilon_1,\ldots,
\varepsilon_n$. One may consider a new game postulating that these
controls are independent. Such game will be called the ordinary differential
game associated with the $\varepsilon$-representation of the interactive
game.

Let us consider now an arbitrary ordinary differential game with the 
evolution equations
$$\dot\varphi=\Phi(\varphi,u_1,u_2,\ldots,u_n),$$
where $\varphi$ characterizes the state of the system and $u_i$ are
the ordinary controls. Let us fix a player. For simplicity of notations we
shall suppose that it is the first one. As a rule the players have their
algorithms of predictions of the behaviour of other players. For a fixed
moment $t_0$ of time let us consider the prediction of the first player for
the game. It consists of the predicted controls $u^\circ_{[t_0];i}(t)$ 
($t>t_0$; $i\ge2$) of all players and the predicted evolution of the system 
$\varphi^\circ_{[t_0]}(t)$. Let us fix $\Delta t$ and consider the real and 
predicted controls for the moment $t_0+\Delta t$. Of course, they may be
different because other players use another algorithms for the game prediction.
One may interpret the real controls $u_i(t)$ ($t=t_0+\Delta t$; $i\ge2$) of 
other players as interactive ones whereas the predicted controls 
$u^\circ_{[t_0];i}(t)$ as pure ones, i.e. to postulate their relation in the 
form:
$$u_i(t)=u_i(u^\circ_{[t-\Delta t];i}(t);\varphi^\circ_{[t_0]}(\tau)|\tau\le t).$$
In particular, the feedbacks may be either reduced to differential form via
the introducing of the intention fields or simply postulated to be differential.
Thus, we constructed an interactive game from the initial ordinary game.
One may use $\varphi(\tau)$ as well as $\varphi^\circ_{[t_0]}(\tau)$ 
in the feedbacks. 

Note that the controls of the first player may be also considered as
interactive if the corrections to the predictions are taken into account
when controls are chosen.
              
The obtained construction may be used in practice to make more adequate
predictions. Namely, {\sl a posteriori\/} analysis of the differential
interactive games allows to make the short-term predictions in such games.
One should use such predictions instead of the initial ones. Note that
at the moment $t_0$ the first player knows the pure controls of other
players at the interval $[t_0,t_0+\Delta t]$ whereas their real freedom
is interpreted as an interactivity of their controls $u_i(t)$. So it is
reasonable to choose $\Delta t$ not greater than the admissible time depth 
of the short-term predictions. 

Na{\"\i}vely, the proposed idea to improve the predictions is to consider
deviations of the real behaviour of players from the predicted ones as
a result of the interactivity, then to make the short-term predictions
taking the interactivity into account and, thus, to receive the 
corrections to the initial predictions. Such corrections may be regarded
as ``psychological'' though really they are a result of different methods
of predictions used by players.

\remark{Remark 3}
The interpretation of the ordinary differential game as an interactive game
also allows to perform the strategical analysis of interactive games. Indeed,
let us consider an arbitrary differential interactive game $A$. Specifying its
$\varepsilon$-representation one is able to construct the associated 
ordinary differential game $B$ with the doubled number of players. Making
some predictions in the game $B$ one transform it back into an interactive 
game $C$. Combination of the strategical long-term predictions in the game 
$B$ with the short-term predictions in $C$ is often sufficient to obtain the 
adequate strategical prognosis for $A$.
\endremark

\remark{Remark 4} The interpretation of the ordinary differential game as
an interactive game is especially useful in situations, when the goals of
players are not known precisely to each other and some more or less rough
suppositions are made.
\endremark

Let's expose another way to obtain an interactive game from the ordinary
differential game. If one has an ordinary differential game with the
states $\varphi(t)$ and the controls $u_i(t)$ it is possible to introduce
some filtering procedure for controls, for instance, the exclusion of the
high frequency components or the separation of frequencies from some special
set. The result will be denote by $u^\circ_i(t)$. One may suspect that 
$u_i^\circ(t)=u_i^\circ(u_i(t);\varphi(t),\dot\varphi(t),\ddot\varphi(t),
\ldots,\varphi^{(k)}(t))$ and use {\sl a posteriori\/} analysis for 
estimation of the feedback. In general, the filtering procedure may depend 
on $\varphi(t)$ and its derivatives but in any case the interactive game is 
constructed. Also one may consider the integrodifferential form of feedbacks:
$u_i^\circ(t)=u_i^\circ(u_i(t),\left.[\varphi(\tau)]\right|_{\tau\leqslant 
t})$.

Note that both versions to unravel an interactivity of the ordinary
differential game may be combined.

\remark{Remark 5} The unraveling of the interactivity of ordinary 
differential games is an essential step for the understanding of 
the complex psychophysiological processes. The next step is to describe
if necessary the interactivity in terms of intention fields and to unravel 
their structure (to solve the inverse problem of representation theory). It 
may be useful to interpret the intention fields as fields of interactive
forces. This scheme may be applied to various interesting psychophysiological 
phenomena usually exposed in terms of Eastern (presumably, Chineese and 
Indian) traditions. Note that their use in medicine is certainly a tactical 
problem (see below).
\endremark

\subhead 1.2. Dialogues and verbalizable games [3]\endsubhead
Let us now expose the interactive game formalism for a description of
dialogues as psycholinguistic phenomena [3]. First of all, note that one is
able to consider interactive games of discrete time as well as interactive
games of continuous time above.

\definition{Defintion 3A (the na{\"\i}ve definition of dialogues) [3]}
The {\it dialogue\/} is a 2-person interactive game of discrete time with 
intention fields of continuous time.
\enddefinition

The states and the controls of a dialogue correspond to the speech whereas 
the intention fields describe the understanding. 

Let us give the formal mathematical definition of dialogues now.

\definition{Definition 3B (the formal definition of dialogues) [3]}
The {\it dialogue\/} is a 2-person interactive game of discrete time of 
the form
$$\varphi_n=\Phi(\varphi_{n-1},\vec v_n,\xi(\tau)| 
t_{n-1}\!\leqslant\!\tau\!\leqslant\!t_n).\tag4$$
Here $\varphi_n\!=\!\varphi(t_n)$ are the states of the system at the
moments $t_n$ ($t_0\!<\!t_1\!<\!t_2\!<\!\ldots\!<\!t_n\!<\!\ldots$), 
$\vec v_n\!=\!\vec v(t_n)\!=\!(v_1(t_n),v_2(t_n))$ are the interactive 
controls at the same moments; $\xi(\tau)$ are the intention fields of 
continuous time with evolution equations
$$\dot\xi(t)=\Xi(\xi(t),\vec u(t)),\tag5$$
where $\vec u(t)=(u_1(t),u_2(t))$ are continuous interactive controls with 
$\varepsilon$--represented couplings of feedbacks:
$$u_i(t)=u_i(u_i^\circ(t),\xi(t);\varepsilon_i(t)).$$
The states $\varphi_n$ and the interactive controls $\vec v_n$ are certain
{\sl known\/} functions of the form
$$\aligned
\varphi_n=&\varphi_n(\vec\varepsilon(\tau),\xi(\tau)| 
t_{n-1}\!\leqslant\!\tau\!\leqslant\!t_n),\\
\vec v_n=&\vec v_n(\vec u^\circ(\tau),\xi(\tau)|
t_{n-1}\!\leqslant\!\tau\!\leqslant\!t_n).
\endaligned\tag6
$$
\enddefinition

Note that the most nontrivial part of mathematical formalization of dialogues
is the claim that the states of the dialogue (which describe a speech) are 
certain ``mean values'' of the $\varepsilon$--parameters of the intention
fields (which describe the understanding).

\remark{Important}
The definition of dialogue may be generalized on arbitrary number of players
and below we shall consider any number $n$ of them, e.g. $n=1$ or $n=3$, 
though it slightly contradicts to the common meaning of the word ``dialogue''.
\endremark

An embedding of dialogues into the interactive game theoretical picture
generates the reciprocal problem: how to interpret an arbitrary differential
interactive game as a dialogue. Such interpretation will be called the
{\it verbalization}.

\definition{Definition 4 [3]} A differential interactive game of the form
$$\dot\varphi(t)=\Phi(\varphi(t),\vec u(t))$$
with $\varepsilon$--represented couplings of feedbacks 
$$u_i(t)=u_i(u^\circ_i(t),\varphi(t),\dot\varphi(t),\ddot\varphi(t),\ldots
\varphi^{(k)}(t);\varepsilon_i(t))$$
is called {\it verbalizable\/} if there exist {\sl a posteriori\/}
partition $t_0\!<\!t_1\!<\!t_2\!<\!\ldots\!<\!t_n\!<\!\ldots$ and the
integrodifferential functionals
$$\aligned
\omega_n&(\vec\varepsilon(\tau),\varphi(\tau)|
t_{n-1}\!\leqslant\!\tau\!\leqslant\!t_n),\\
\vec v_n&(\vec u^\circ(\tau),\varphi(\tau)|
t_{n-1}\!\leqslant\!\tau\!\leqslant\!t_n)
\endaligned\tag7$$ 
such that
$$\omega_n=\Omega(\omega_{n-1},v_n;\varphi(\tau)|
t_{n-1}\!\leqslant\!\tau\!\leqslant\!t_n).
\tag 8$$
\enddefinition

The verbalizable differential interactive games realize a dialogue in sense
of Def.3.

\remark{Remark 6} One may include $\omega_n$ explicitely into the evolution
equations for $\varphi$
$$\dot\varphi(\tau)=\Phi(\varphi(\tau),\vec u(\tau);\omega_n),\quad 
\tau\in[t_n,t_{n+1}].$$
as well as into the feedbacks and their couplings.
\endremark

The main heuristic hypothesis is that all differential interactive games
``which appear in practice'' are verbalizable. The verbalization means that 
the states of a differential interactive game are interpreted as intention 
fields of a hidden dialogue and the problem is to describe such dialogue 
completely. If a differential interactive game is verbalizable one 
is able to consider many linguistic (e.g. the formal grammar of a related 
hidden dialogue) or psycholinguistic (e.g. the dynamical correlation of 
various implications) aspects of it.

During the verbalization it is a problem to determine the moments $t_i$. A 
way to the solution lies in the structure of $\varepsilon$-representation.
Let the space $E$ of all admissible values of $\varepsilon$-parameters be
a CW-complex. Then $t_i$ are just the moments of transition of the 
$\varepsilon$-parameters to a new cell. 

\subhead 1.3. Tactical games [4]\endsubhead
Tactics as it will be defined below is derived from two independent concepts:
the parametric interactive games and external controls on one hand and
the comments to dialogues on another hand.

First of all, note that an interactive game may depend on the additional 
paramaters. Such dependence is of two forms. First, parameters may appear in 
the evolution equations:
$$\dot\varphi=\Phi(\varphi,u_1,u_2,\ldots,u_n;\lambda).\tag9A$$
Here, $\lambda$ is a collective notation for parameters. Second, parameters
may appear in feedbacks:
$$u_i(t)=u_i(u_i^\circ(t),\varphi(t),\ldots,\varphi^{(k)}(t);\lambda).\tag9B$$
The dependence of $u_i$ on $\lambda$ is either unknown (incompletely known) or 
known. The least means that $\partial u_i/\partial\lambda$ may be expressed via 
$u_i$ as a function of other variables (such expression are integrodifferential 
on these variables). Both variants of parametric dependence of interactive game 
may be combined together.

The additional paramaters may realize the external controls. In this
situation they depend on time:
$$\lambda=\lambda(t).$$
In practice, such situation appear in the teaching systems. The parameters
are interpreted as controls of a teacher. This example is rather typical.
It shows that the controls $\lambda(t)$ may be considered as ``slow''
whereas the interactive controls $u_i(t)$ as ``quick''.

Of course, one is able to introduce the slow controls $\lambda(t)$, which
belong to the same players as the interactive controls $u_i(t)$ or to their
coalitions. And, certainly, the slow controls of discrete time may be
considered. One may suspect that the discrete time controls $\lambda_n$
realize a convenient approximation for the slow controls $\lambda(t)$, which
is timer in practice.

The slow controls may be interactive.

If dependence of $u_i$ on $\lambda$ is known and one consider the 
$\varepsilon$-representation of feedbacks it is either postulated that 
$\varepsilon$-parameters do not depend on $\lambda$ or claimed that 
$\partial\varepsilon/\partial\lambda$ is expressed via $\varepsilon$ as 
a function of other arguments.
                                                   
Now let us define comments to the dialogue.
Let 
$$\omega_n=\Omega(\omega_{n-1},\vec v_n,\xi(\tau)|
t_{n-1}\!\leqslant\!\tau\!\leqslant\!t_n)$$
be the $n$-person dialogue with the discrete time interactive controls 
$\vec v_n$ and the intention fields governed by the evolution equations
$$\dot\xi(t)=\Xi(\xi(t),\vec u(t)),$$
where $\vec u(t)$ are the continuous interactive controls with
$\varepsilon$--represented couplings of feedbacks:
$$u_i(t)=u_i(u_i^\circ(t),\xi(t);\varepsilon_i(t)).$$
The states $\varphi_n$ and the interactive controls $\vec v_n$ are
expressed as
$$\aligned
\omega_n=&\omega_n(\vec\varepsilon(\tau),\xi(\tau)|
t_{n-1}\!\leqslant\!\tau\!\leqslant\!t_n),\\
\vec v_n=&\vec v_n(\vec u^\circ(\tau),\xi(\tau)|
t_{n-1}\!\leqslant\!\tau\!\leqslant\!t_n).
\endaligned
$$

The discrete time {\it comments\/} $\vartheta_n$ to
the dialogue are defined recurrently as
$$\vartheta_n=\Theta(\vartheta_{n-1},\omega_n,v_n).\tag10$$

Comments to the dialogue at the fixed moment $t_n$ contain various
information on the dialogue. For instance, one may to raise a problem
to restore some features of a dialogue from certain comments or alternatively
what features of a dialogue may be restored from the fixed comment.

The main difference of the comments $\vartheta_n$ from the states $\omega_n$
is the absence of expressions of the first via $\vec\varepsilon(\tau)$ and
$\xi(\tau)$ ($t_{n-1}\!\leqslant\!\tau\!\leqslant\!t_n$).

Comments are applied to the verbalizable games in the same way.
                                                               
Tactical games combine mechanisms of parametric interactive games and 
comments to the dialogue (verbalizable game).

\definition{Definition 5} The {\it tactical game\/} is a parametric
verbalizable game with comments, in which the parameters are of discrete
time and coincide with the comments.
\enddefinition

It is really wonderful that such simple definition is applicable to a very
huge class of phenomena. However, it is so! As it was marked above virtually
all known forms of human activity such as scientific researches and
economics, sport or military actions, medicine, fine or martial arts, 
literature and music, theatre and dance, psychotherapy and even magic may be 
regarded as certain tactical games. Trying to improve the model I have no 
found any wider concept, whose using is necessary and effective, whereas 
the notion of tactical game may describe these phenomena very correctly.

\remark{Remark 7} The pairs $(v_n,\vartheta_n)$ of discrete time interactive
controls and the comments will be called the {\it tactical actions}, whereas
the continuous time interactive controls $u_i(t)$ will be called the {\it 
instant actions}. The tactical actions may be involved as in the evolution 
equations as in the interactivity.
\endremark

Now we shall describe some operations over the tactical games (the tactical
interaction, the tactical synthesis and the tactical extension).

Let us consider two tactical games defined by the evolution equations
$$\dot\varphi_1=\Phi_1(\varphi_1,\vec u_1;\vartheta_1)$$
and 
$$\dot\varphi_2=\Phi_2(\varphi_2,\vec u_2;\vartheta_2)$$
with $\varepsilon$--represented couplings of feedbacks
$$u_{1,i}=u_{1,i}(u^\circ_{1,i},\varphi_1,\dot\varphi_1,
\ddot\varphi_1,\ldots \varphi^{(k)}_1;\varepsilon_{1,i},
\vartheta_1)$$
and
$$u_{2,i}=u_{2,i}(u^\circ_{2,i},\varphi_2,\dot\varphi_2,
\ddot\varphi_2,\ldots \varphi^{(k)}_2;\varepsilon_{2,i},
\vartheta_2).$$
The integrodifferential functionals (7) have the form
$$\aligned
\omega_{j,n}&(\vec\varepsilon_j(\tau),\varphi_j(\tau)|
t_{n-1}\!\leqslant\!\tau\!\leqslant\!t_n),\\
\vec v_{j,n}&(\vec u_j^\circ(\tau),\varphi_j(\tau)|
t_{n-1}\!\leqslant\!\tau\!\leqslant\!t_n)
\endaligned$$
and the relations (8)
$$\omega_{j,n}=\Omega_j(\omega_{j,n-1},v_{j,n};\varphi_j(\tau)|
t_{n-1}\!\leqslant\!\tau\!\leqslant\!t_n)$$
hold ($j=1,2$). The comments $\vartheta_1$ and $\vartheta_2$ are defined
recurrently as
$$\vartheta_{1,n}=\Theta_1(\vartheta_{1,n-1},\omega_{1,n},v_{1,n})$$
and
$$\vartheta_{2,n}=\Theta_2(\vartheta_{2,n-1},\omega_{2,n},v_{2,n}).$$

\define\tni{\operatorname{int}}
The {\it tactical interaction\/} is realized by the addition of the
interaction terms into the recurrent formulas for $\vartheta_j$ to produce
the interdetermination of
comments:
$$\vartheta_{1,n}=\Theta_1(\vartheta_{1,n-1},\omega_{1,n},v_{1,n})+
\tilde\Theta_{1,2}^{\tni}(\vartheta_{1,n-1},\vartheta_{2,n-1},\omega_{1,n},
v_{1,n})\tag11A$$
and
$$\vartheta_{2,n}=\Theta_2(\vartheta_{2,n-1},\omega_{2,n},v_{2,n})+
\tilde\Theta_{2,1}^{\tni}(\vartheta_{2,n-1},\vartheta_{1,n-1},\omega_{2,n},
v_{2,n}).\tag11B$$
         
Let us consider $N$ control systems represented as tactical games defined
by the evolution equations
$$\dot\varphi_j=\Phi_j(\varphi_j,\vec u_j;\vartheta_j)$$
($j=1,2,\ldots N$) with $\varepsilon$--represented couplings of feedbacks
$$u_{j,i}=u_{j,i}(u^\circ_{j,i},\varphi_j,\dot\varphi_j,
\ddot\varphi_j,\ldots \varphi^{(k)}_j;\varepsilon_{j,i},
\vartheta_j).$$
The integrodifferential functionals (7) have the form
$$\aligned
\omega_{j,n}&(\vec\varepsilon_j(\tau),\varphi_j(\tau)|
t_{n-1}\!\leqslant\!\tau\!\leqslant\!t_n),\\
\vec v_{j,n}&(\vec u_j^\circ(\tau),\varphi_j(\tau)|
t_{n-1}\!\leqslant\!\tau\!\leqslant\!t_n)
\endaligned$$
and the relations (8)
$$\omega_{j,n}=\Omega_j(\omega_{j,n-1},v_{j,n};\varphi_j(\tau)|
t_{n-1}\!\leqslant\!\tau\!\leqslant\!t_n)$$
hold. The comments $\vartheta_j$ are defined recurrently as
$$\vartheta_{j,n}=\Theta_j(\vartheta_{j,n-1},\omega_{j,n},v_{j,n}).\tag12$$

The {\it tactical synthesis\/} is realized by the redefinition of the
recurrent formulas for $\vartheta_j$ to produce the unification of
comments:
$$\vartheta_{j,n}=\tilde\Theta_j(\vartheta_{1,n-1},\ldots\vartheta_{N,n-1},
\omega_{1,n},\ldots,\omega_{N,n},v_{1,n},\ldots v_{N,n}).$$
The functions $\tilde{\boldsymbol\Theta}\!=\!(\tilde\Theta_1,\ldots\tilde
\Theta_N)$ determines the synthesis.It may has various internal structure,
which is characterized by the set of real arguments of functions 
$\tilde\Theta_j$ and their hierarchical structure. It presupposes that 
$\tilde\Theta_j$ depend not on all triples $(\vartheta_j,\omega_j,v_j)$ and 
various triples may appear in $\tilde\Theta_j$ in coalitions of different form 
and nature. One may think that the functions $\tilde{\boldsymbol\Theta}\!=
\!(\tilde\Theta_1,\ldots\tilde\Theta_N)$ are constructed from the functions
$\boldsymbol\Theta=(\Theta_1,\ldots\Theta_N)$ using some operations, which
realize the synthesis. Tactical interaction is a form of tactical synthesis
of two games.

If under the tactical synthesis of two games (in particular, it may be
a tactical interaction) $\tilde\Theta_1=\Theta_1$ we shall say that
the {\it tactical extension\/} of the first game is realized.

\head II. Dialectics and dialectical games\endhead

To unravel the abstract foundations of tactics is much more than to define
it. Such foundations should provide the description of tactical phenomena
that presupposes a possibility of their repreating and reproducing, whereas
definitions only established the rules of games. The following paragraph is
an attempt to initiate the investigation of such foundations. The crucial
role is played by the game theoretical definition of dialectics.

\definition{Definition 6} {\it Dialectics\/} is a logical self-describing
tactical game. 
\enddefinition

"Logical" means that the game is purely cognitive and does not need obligatory
any material substrate or, otherwise, that its objects are notions. 
"Self-describing" means that all notions as well as procedures of their
transformation, which constitute the game evolution, are derived from 
themselves and that self-description as a definition of all appearing objects
from themselves is simultaneously the goal of the game.

The main question is whether dialectics exists. Its existence should be 
considered as a postulate.

\proclaim{Postulate} Dialectics exists.
\endproclaim

Apparently this postulate can not be derived in any way from the standard
mathematical foundations such as set theory, mathematical logic and
arithmetics. Moreover, there are no reasons to believe that it is possible
to prove formally (i.e. by means of mathematical logic) that it does not 
contradict to them. However, if the postulate is adopted its compatibility
with mathematical logic may be received dialectically.

Another question is one of the equivalence of various forms of dialectics.
It is not a fact that dialectics is unique.

Now let us define dialectical games.

\definition{Definition 7} A {\it dialectical game\/} is a tactical extension
of dialectics.
\enddefinition

It means that dialectical games use some dialectical procedures combined
with others, otherwords, have elements of self-description. Dialectical
games are not obligatory logical. For instance, one may consider dialectical 
games of perception, dialectical computer games, etc.

Note that dialectics as a tactical game itself is not given {\sl before\/}
but is self-constructing {\sl during\/} the dialectical game. It is its
difference from the formal logic.

\remark{Remark 8} Dialectics allows to give self-consistent descriptions of 
any tactical games and such descriptions are dialectical games. So one may
think dialectics as a formal cause of descriptibility of tactical phenomena.
\endremark

Thus, if tactics is thinking as "an art to manipulate the unknown, which is
manifested by the interactivity, without making it known" [4], dialectics 
should be regarded as "an art to comprehend these manipulations" and, 
moreover, as "a kinaesthetic art to comprehend the unknown itself via 
manipulations over it", dialectics is an abstract and kinaesthetic 
contemplation of the unknown. Hence, the dialetical games may be an effective 
tools for analysis and controlling of the internal dialogue and the reflection
because dialectics may be regarded as an {\sl interiorization of tactics}.

In fact, tactical games tends to be dialectical because any purposeful
activity presuppose some form of its understanding based on its description.
Because description of a tactical phenomena is also a tactical game, to make 
such description self-consistent one should use dialectics (or already 
constructed dialectical game). Hence, tactics is effectiveness as much as
it uses dialectics.

Let us describe a general form of the dialectical game. Such game is defined 
by the evolution equations
$$\dot\varphi=\Phi(\varphi,\vec u;\vartheta)$$
with $\varepsilon$--represented couplings of feedbacks
$$u_i=u_i(u^\circ_i,\varphi,\dot\varphi,\ddot\varphi,\ldots \varphi^{(k)};
\varepsilon_i,\vartheta)$$
The integrodifferential functionals (7) have the form
$$\aligned
\omega_n&(\vec\varepsilon(\tau),\varphi(\tau)|
t_{n-1}\!\leqslant\!\tau\!\leqslant\!t_n),\\
\vec v_n&(\vec u^\circ(\tau),\varphi(\tau)|
t_{n-1}\!\leqslant\!\tau\!\leqslant\!t_n)
\endaligned$$
and the relations (8)
$$\omega_n=\Omega(\omega_{n-1},v_n;\varphi(\tau)|
t_{n-1}\!\leqslant\!\tau\!\leqslant\!t_n).$$
The comments $\vartheta_n$ are defined recurrently as
$$\vartheta_n=\Theta(\vartheta_{n-1},\delta_n,\omega_n,v_n),\tag13$$
where $\{\delta_n\}$ are self-describing objects of dialectics.

\remark{Remark 9} In the dialectical game the tactical actions (the pairs
$(v_n,\vartheta_n)$) are constructed using the dialectical objects $\delta_n$
besides $v_n$ and $\omega_n$ so dialectics appears as a logic of tactics.
\endremark

\remark{Remark 10} Often a dialectical game is a tactical extension of
both dialectics and a concrete tactical game of any nature.
\endremark

\remark\nofrills{Remark 11 (for a discussion):}\ \eightpoint Dialectics may be 
regarded not only as a formal cause of descriptibility of tactical phenomena 
but also as a formal cause of these phenomena themselves. For an active cause 
of tactical (interactive) phenomena the notion "m\B ay\B a" was adopted. One 
may think that both causes form the same {\sl linguodynamical or logodynamical
cause\/} of tactics, and for such linguodynamical unity of dialectics as 
a logical self-descripting tactical game and m\B ay\B a as a universal 
dynamical principle of all subject-object relations the related term 
{\sl "lil\B a"\/} is very convenient. Because m\B ay\B a is not neither
object nor subject, and it is not active, there are no any energetical or 
dialogical aspects of lil\B a itself, besides ones between real subjects and
objects. Lil\B a is thought as purely dynamical, kinaesthetical and 
contemplative dialectics of m\B ay\B a. Note that m\B ay\B a is objectivized 
in various forms: as intention fields or as fields of interactive forces. 
All these objects may be visualized. The visualizations of lil\B a are 
the dialectical perception games, combining visual, dynamical, kinaesthetical 
and dialectically contemplative aspects.
\endremark
\tenpoint

\head III. Dialectical games and representation theory\endhead

Let us describe a type of dialectical games, which may be considered
simultaneously as illustrations of general concepts and as models
for the internal dialogue and the reflection. 

First of all, some additional constructions, which realize a link between
control and representation theory, are necessary.
\newpage

\subhead 3.1. Representative dynamics [5]\endsubhead

\definition{Definition 8} Let $\bold X=\bold X(t)=(X_1(t),\ldots X_m(t))$
($X_i(t)\in\Mat_n(\Bbb C)$) be the time-dependent vector of $m$ complex 
$n\times n$ matrices. The {\it representative dynamics\/} is a controlled
system (with constraints) of the form
$$\dot\bold X(t)=F(\bold X(t),a(t))\tag14$$
with the {\sl fixed\/} initial data $\bold X(t_0)$, where the control 
parameter $a(t)=(\frak A(t),\bold e(t))$ is the pair of {\sl any\/} 
associative algebra $\frak A(t)$ from the fixed class of such algebras 
$\Bbb A$, $\bold e(t)=(e_1(t),\ldots e_m(t))$ is {\sl any\/} set of algebraic 
generators of the algebra $\frak A(t)$ (one may claim $\bold e(t)$ to be an 
algebraic basis in $\frak A(t)$) such that the mapping $e_i(t)\mapsto X_i(t)$ 
may be extended to the representation $T(t):\frak A(t)\mapsto\Mat_n(\Bbb C)$ 
of the algebra $\frak A(t)$ in the matrix algebra $\Mat_n(\Bbb C)$ (this is a 
constraint on the control $a(t)$).
\enddefinition

Certainly, the claim that (1) is a representative dynamics restricts the
choice of the function $F$ and initial data $\bold X(t_0)$ because for
each moment $t$ {\sl any\/} admissible choice of the pair $(\frak A(t),
\bold e(t))$ should provide that the set of admissible pairs will not be 
empty in future. 

\remark{Remark 12} It is possible to consider the infinite dimensional 
algebras of operators instead of $\Mat_n(\Bbb C)$ or their subalgebras
(e.g. commutative algebras of functions if all algebras $\frak A(t)$ are
commutative).
\endremark

\remark{Remark 13} Let us consider the following equivalence on the set 
$\Cal A$ of all admissible $a=(\frak A,\bold e)$. The pairs $a_1=(\frak A_1,
\bold e_1)$ and $a_2=(\frak A_2,\bold e_2)$ will be equivalent iff the 
algebras $A_1$ and $A_2$ are isomorphic under an isomorphism which maps the 
linear space $V_1$ spanned by the elements of $\bold e_1$ onto the linear 
space $V_2$ spanned by the elements of $\bold e_2$. Then the equivalence
divides the time interval $[t_0,t_1]$, on which the representative dynamics is 
considered, onto the subsets, on which the pairs $a(t)$ are equivalent.
\endremark

\remark{Remark 14} Representative dynamics combines structural and functional
features. The first are accumulated in the class $\Bbb A$ and the least are
expressed by the function $F$. Both are interrelated. The situation is similar
to one in the group theory of special functions. However, the difference is 
essential: in the representative dynamics the functional aspects are not 
derived from the structural ones and hve an independent origin.
\endremark

Let us now describe the dynamical inverse problem of representation theory for
controlled systems following the general ideology of the inverse problems
of representation theory [6].

\definition{Definition 9} Let 
$$\dot x=\varphi(x,u),\tag15$$
be the controlled system, where $x$ is the time-dependent $m$-dimensional
complex vector and $u$ is the control parameter. {\it Dynamical inverse 
problem of representation theory for the controlled system\/} (15) is to 
construct a representative dynamics
$$\dot\bold X=F(\bold X,a)$$
and the function
$$a=a(u,x)\quad\text{such that}\quad \varphi(x,u)=f(x,a(u,x)),$$
where the operator function $F$ is defined by the Weyl (symmetric) symbol $f$ 
as a function of $m$ non-commuting variables $X_1,\ldots X_m$.
\enddefinition

\remark{Remark 15} If the controls are absent and the pair $a(t)=(\frak A(t),
\bold e(t))$ is time-independent, Def.9 is reduced to the definition of the 
dynamical inverse problem of representation theory [6].
\endremark

\remark{Remark 16} One may consider dynamical inverse problem of representation
theory for games, the interactively controlled systems and interactive games.
\endremark

\remark{Remark 17} If the function $\varphi$ contained some constants
$c_\alpha\in\Bbb C$ then one may interpret them as time-independent variables
and include the matrices $C_\alpha\in\Mat_n(\Bbb C)$ instead of them in the
operator function $F$ (compare with the quantization of constants [6]).
\endremark

\subhead 3.2. Tactical representative dynamics as dialectical game\endsubhead
Note that a representative dynamics may be considered as a parametric one.
The parameter is a class $\Bbb A$ of associative algebras. Such view allows
to define the tactical representative dynamics.

\definition{Definition 10} The {\it tactical representative dynamics\/}
is a tactical game, which is a representative dynamics as a verbalizable
interactive game, with the variable classes $\Bbb A_n$ of associative
algebras as the comments.
\enddefinition

The change of the class $\Bbb A_n$ may be caused by the fact that
the equations of representative dynamics become insolvable in the class
$\Bbb A_{n-1}$.

One may describe solutions of many mathematical problems as tactical 
representative dynamics. The keypoint is an algorithm of the constructing
of comments, i.e. how to derive a class $\Bbb A_n$ from the preceeding one
$\Bbb A_{n-1}$. Though in the simplest case the recurrent formulas for 
comments are represented as
$$\vartheta_n=\Theta(\vartheta_{n-1}),$$
where $\vartheta_n$ is the pair $(\Bbb A_n,\eta_n)$ of the class $\Bbb A_n$
and additional parameters $\eta_n$, the general form is 
$$\vartheta_n=\Theta(\vartheta_{n-1},\omega_n,v_n)$$
as it is in the definition of a tactical game or explicitely
$$\left\{\aligned
\Bbb A_n=&\bold A(\Bbb A_{n-1},\eta_{n-1},\omega_n,v_n)\\
\eta_n=&\bold H(\eta_{n-1},\Bbb A_{n-1},\omega_n,v_n)
\endaligned\right.$$
However, it is intuitively clear that the dialectical objects should be 
involved into this procedure because there are no any self-consistent
ways to obtain any non-trivially and principally new class $\Bbb A_n$ from 
the preceeding one $\Bbb A_{n-1}$ so $\Bbb A_n$ should be self-derived from
$\Bbb A_{n-1}$ and the transition $\Bbb A_{n-1}\Rightarrow\Bbb A_n$ should
be described in the internal terms of the class $\Bbb A_n$. So the tactical 
representative dynamics corresponding to the solution of a mathematical 
problem should be a dialectical game and the algorithm of the comments' 
constructing should have the form
$$\vartheta_n=\Theta(\vartheta_{n-1},\delta_n,\omega_n,v_n),$$
where $\delta_n$ are dialectical objects. Explicitely,
$$\left\{\aligned
\Bbb A_n=&\bold A(\Bbb A_{n-1},\eta_{n-1},\delta_n,\omega_n,v_n)\\
\eta_n=&\bold H(\eta_{n-1},\Bbb A_{n-1},\delta_n,\omega_n,v_n)
\endaligned\right.$$
Such tactical representative dynamics will be called {\it dialectical}.
The expressions for $\Bbb A_n$ describe its synthesis, which is the
dialectical synthesis of a class $\Bbb A_n$ of associative algebras.

\remark{Remark 18} The dialectical representative dynamics may be considered
as a model for the internal dialogue and reflection of a human. Here, the
representations in mathematical and psychological sense coincide. The
associative algebras $\frak A(t)$ from the classes $\Bbb A_n$ describe the 
cognitive aspect of the process. The recurrent formulas for comments 
correspond to the deduction. 
\endremark

\remark{Remark 19} The tactical (dialectical) representative dynamics 
solves the inverse problem of representation theory for the tactical
(dialectical) games.
\endremark

\remark{Remark 20} One can describe a dialectical object $\delta$ by
the concrete transitions (syntheses) $\Bbb A_n\Rightarrow\Bbb A_{n+1}$, 
which it governs. 
\endremark

\head Conclusions\endhead

Thus, the tactical game theoretical interpretation of dialectics is
given. Dialectical games are considered as abstractly as well as models
of the internal dialogue and reflection. The models related to the
representation theory (representative dynamics) are specially investigated
in detail, they correlate with the hypothesis on the dialectical features
of human thinking in general and mathematical thought (the constructing of
a solution of mathematical problem) in particular. 

Also the discussed subjects are significant for the understanding of tactical 
features of various psychophysiological processes and their couplings with 
cognitive ones during the functioning of the videocognitive and analogous 
integrated interactive systems and, thus in fact, for the clarification of 
the psychophysical nature of cognitive activity.

\Refs
\roster
\item"[1]" Isaaks R., Differential games. Wiley, New York, 1965;\newline
Owen G., Game theory, Saunders, Philadelphia, 1968;\newline
Vorob'ev N.N., Current state of the game theory, Uspekhi Matem. Nauk 15(2) 
(1970) 80-140 [in Russian].
\item"[2]" Juriev D., Interactive games and representation theory. I,II.
E-prints: math.FA/9803020, math.RT/9808098;\newline
Juriev D., Games, predictions, interactivity. E-print: math.OC/9906107.
\item"[3]" Juriev D., Interactive games, dialogues and the verbalization.
E-print: math.OC/9903001.
\item"[4]" Juriev D., Tactical games \& behavioral self-organization.
E-print: math.HO/9911120.
\item"[5]" Juriev D., Representative dynamics. E-print: math.OC/9911147.
\item"[6]" Juriev D.V., An excursus into the inverse problem of representation
theory [in Russian]. Report RCMPI-95/04 (1995) [e-version: mp\_arc/96-477].
\endroster
\endRefs
\enddocument